\documentclass[12pt]{amsart}
\usepackage{amssymb} 

\newcommand{\uppermu}{{\overline{\mu}}}
\newcommand{\lowermu}{{\underline{\mu}}}
\newcommand{\weak}{{\operatorname{weak}}}

\newcommand{\BOX}{{\operatorname{Box}}}

\newcommand{\Fqt}{{\F_q[t]}}
\newcommand{\re}{\operatorname{Re}}

\newcommand{\twoline}[2]{%
  {\mathchoice{\genfrac{}{}{0pt}0{#1}{#2}}%
              {\genfrac{}{}{0pt}0{#1}{#2}}%
              {\genfrac{}{}{0pt}1{#1}{#2}}%
              {\genfrac{}{}{0pt}2{#1}{#2}}}}

\newcommand{\Q}{{\mathbf Q}}

\newcommand{\Z}{{\mathbf Z}}

\newcommand{\Aff}{{\mathbf A}}

\newcommand{\calQ}{{\mathcal Q}}

\newcommand{\calR}{{\mathcal R}}
\newcommand{\calS}{{\mathcal S}}
\newcommand{\calT}{{\mathcal T}}
\newcommand{\calZ}{{\mathcal Z}}

\newcommand{\PP}{{\mathbf P}}

\newcommand{\OO}{{\mathcal O}}

\newcommand{\fraks}{{\mathfrak s}}
\newcommand{\pp}{{\mathfrak p}}

\newcommand{\F}{{\mathbf F}}

\newcommand{\Aut}{\operatorname{Aut}}

\newcommand{\Div}{\operatorname{Div}}

\newcommand{\Spec}{\operatorname{Spec}}

\newcommand{\isom}{\simeq}
\newcommand{\del}{\partial}

\newtheorem{theorem}{Theorem}[section]
\newtheorem{lemma}[theorem]{Lemma}
\newtheorem{cor}[theorem]{Corollary}

\theoremstyle{definition}

\newtheorem{question}[theorem]{Question}

\theoremstyle{remark}
\newtheorem{rem}{Remark$\!\!$}

\begin{document}

\title{Squarefree values of multivariable polynomials}
\subjclass[2000]{Primary 11C08}
\author{Bjorn Poonen}
\thanks{This research was supported by NSF grant DMS-9801104
          and a Packard Fellowship.  Part of the research was done
	while the author was enjoying the hospitality of the
	Universit\'e de Paris-Sud.  The final version of this article
	will be published in the {\em Duke Mathematical Journal},
	published by Duke University Press.}
\address{Department of Mathematics, University of California, Berkeley, CA 94720-3840, USA}
\email{poonen@math.berkeley.edu}
\date{September 16, 2002}

\begin{abstract}
Given $f \in \Z[x_1,\dots,x_n]$,
we compute the density of $x \in \Z^n$ such that $f(x)$ is squarefree,
assuming the $abc$ conjecture.
Given $f,g \in \Z[x_1,\dots,x_n]$,
we compute unconditionally the density of $x \in \Z^n$
such that $\gcd(f(x),g(x))=1$.
Function field analogues of both results are proved unconditionally.
Finally, assuming the $abc$ conjecture, given $f \in \Z[x]$,
we estimate the size of the image of $f(\{1,2,\dots,n\})$
in $\left(\Q^*/\Q^{*2} \right) \cup \{0\}$.
\end{abstract}

\maketitle

\section{Introduction}
\label{introduction}

An integer $n$ is called {\em squarefree} if for all prime numbers $p$
we have $p^2 \nmid n$ (that is, $p^2$ does not divide $n$).
Heuristically, one expects that if one chooses a positive integer $n$
``at random,'' then for each prime $p$, 
the ``probability'' that $p^2 \nmid n$ equals $1-p^{-2}$;
and the assumption that these probabilities are ``independent''
leads to the guess that the density of squarefree positive integers
equals
	$$\prod_{\text{prime } p} (1 - p^{-2}) = \zeta(2)^{-1} = 6/\pi^2,$$
where $\zeta(s)$ is the Riemann zeta function, defined by
	$$\zeta(s) = \sum_{n \ge 1} n^{-s} 
	= \prod_{\text{prime }p} (1-p^{-s})^{-1}$$
for $\re s>1$.
One can formulate this guess precisely by
defining the density of a set of positive integers $S$ as
	$$\mu(S) := \lim_{B \rightarrow \infty} \frac{\# ( S \cap [1,B])}{B}.$$
In fact, 
the guess can be proved by simple sieve techniques~\cite[\S18.6]{hardy-wright}.

Now suppose that $f(x)$ is a polynomial with integer coefficients,
and let $S$ be the set of positive integers $n$ for which $f(n)$ is squarefree.
This time one guesses that the density of $S$
equals $\prod_{\text{prime }p} (1-c_p/p^2)$
where $c_p$ equals the number of integers $n \in [0,p^2-1]$
for which $p^2 \mid f(n)$.
When $\deg f \le 2$, a simple sieve again shows that the guess is correct.
When $\deg f = 3$, a more complicated argument is needed
(see~\cite{hooley1967}, or, for an improved error term, 
Chapter~4 of~\cite{hooley1976}).
For general $f$ with $\deg f \ge 4$, it is unknown whether the
heuristic conjecture is correct, 
but A. Granville~\cite{granville1998} 
showed that it follows from the $abc$ conjecture.
(Recall that the $abc$ conjecture is the statement that
for any $\epsilon > 0$,
there exists a constant $C=C(\epsilon)>0$ such that if $a,b,c$
are coprime positive integers satisfying $a+b=c$,
then $c < C ({\displaystyle \prod_{p | abc} p} )^{1+\epsilon}$.)
Granville used the $abc$ conjecture in conjunction with Belyi's Theorem,
to bound the number of polynomial values divisible by the square
of a large prime.
He also proved a conditional result for homogeneous polynomials
in two variables, extending some earlier results along these lines,
such as~\cite{greaves1992}.
(See~\cite{granville1998} for more references; some of
these earlier results were unconditional in low degree cases.)

In this paper we generalize Granville's results to 
arbitrary polynomials over $\Z$ in many variables, 
still assuming the $abc$ conjecture.
The proof proceeds by reduction to the one-variable case,
and the $abc$ conjecture is required only because it is used
by Granville; it is not required for the reduction.
Such a fibering argument was used also in~\cite{gouvea-mazur1991}.
One defect of our proof is that it appears not to work
for the most natural generalization of density in the multivariable case:
see Section~\ref{definitions} and the remark following the proof
of Lemma~\ref{squarefreeerror} for more details.
An application of our result is towards estimating,
given a regular quasiprojective scheme $X$ over $\Z$,
what fraction of hypersurface sections of $X$ are regular.
(See~\cite{poonenbertini}.)

If $\F_q$ is a finite field of characteristic $p$,
we prove an analogue for polynomials over $\F_q[t]$ unconditionally,
using a completely different proof,
exploiting the fact that $\F_q[t]$ has an $\F_q[t^p]$-linear derivation.
One application of this result,
suggested by A.~J.~de~Jong~\cite[\S4.22]{dejong2001}, 
is to counting elliptic curves with squarefree discriminant:
see Section~\ref{theorems}.
The case of squarefree values of 
a separable irreducible one-variable polynomial
over $\F_q[t]$ (or more generally $k^{\operatorname{th}}$-power-free
values for polynomials over the ring of regular functions
on any affine curve over $\F_q$)
was proved earlier by K.~Ramsay~\cite{ramsay1992}\footnote{The formula
	for the density in Theorem~1 of~\cite{ramsay1992} should read
	$Z=\prod_{v \not\in S} \left(1 - \rho(kv)/\|v\|^k \right)$.
	The proof there is correct, but the statement is 
	unfortunately misprinted,
	with $\rho(v)$ in place of $\rho(kv)$.}
using a lemma of N.~Elkies involving a derivation.
In Section~\ref{ringsoffunctions},
we sketch a generalization of our result
to multivariable polynomials
over such rings of regular functions.

A related problem asks, given relatively prime polynomials 
$f(x_1,\dots,x_n)$
and $g(x_1,\dots,x_n)$ over $\Z$,
what is the density of $n$-tuples of positive integers
for which the values of $f$ and $g$ are relatively prime?
Again there is a heuristic guess, and it was proved in~\cite{ekedahl1991}
that this guess is correct.
We generalize by using a stronger definition of density
(involving boxes of arbitrary dimensions, instead of only equal
dimensions as considered in~\cite{ekedahl1991})
and by simultaneously proving the function field analogue.
The generalizations are needed
to prove the corresponding results about squarefree values.

Finally, we confirm a guess made in~\cite{granville1998},
namely that for a nonzero polynomial $f(x) \in \Z[x]$,
the size of the image of $\{f(1),f(2),\dots,f(B)\}$
in $\left( \Q^*/\Q^{*2} \right) \cup \{0\}$
is $c_f B + o(B)$ as $B \rightarrow \infty$,
for some constant $c_f$ depending on $f$.
Moreover, we find an explicit formula for $c_f$.
In particular, $c_f=1$ if $f$ is squarefree of degree $\ge 2$.

\section{Definition of density}
\label{definitions}

In Sections \ref{definitions} through~\ref{Fqtsection},
$A$ denotes $\Z$ or $\F_q[t]$ for some prime power $q=p^e$.
Let $K$ denote the fraction field of $A$.
For nonzero $a \in A$ define $|a|:=\#(A/a)$,
and define $|0|=0$.
If $\pp$ is a nonzero prime of $A$, let $|\pp|:=\#(A/\pp)$.
Define
	$$\BOX = \BOX(B_1,\dots,B_n) = 
	\begin{cases}
		\{\, (a_1,\dots,a_n) \in \Z^n : 0<a_i \le B_i 
				\text{ for all $i$} \,\} 
			&\text{if $A=\Z$,} \\
		\{\, (a_1,\dots,a_n) \in A^n : |a_i| \le B_i 
				\text{ for all $i$} \,\} 
			&\text{if $A=\F_q[t]$.} 
	\end{cases}$$
For $\calS \subseteq A^n$,
define
	$$\mu(\calS) := \lim_{B_1,\dots,B_n \rightarrow \infty} 
	\frac{\#(\calS \cap \BOX)}{\# \BOX},$$
and define $\uppermu(\calS)$ and $\lowermu(\calS)$ similarly
using $\limsup$ and $\liminf$ in place of $\lim$.
If a subset $\calS \subseteq \Z^n$ and its $2^n$ reflections
in the coordinate hyperplanes have a common density in this strong sense,
then we can estimate $\#(\calS \cap R)/\#R$ for regions $R$
of many other shapes.
For instance, if $R_B$ is the ball of radius $B$ centered at the origin,
then $\#(\calS \cap R_B)/\#R_B \rightarrow \mu(\calS)$ 
as $B \rightarrow \infty$, since $R_B$ can be approximated by a Boolean
combination of $k$ boxes and their reflections,
with an error of at most $\epsilon_k B^n$ lattice points
for $B$ large relative to $k$,
where $\epsilon_k \rightarrow 0$ as $k \rightarrow \infty$.

In some of our results
we can prove that the density exists only in a weaker sense.
Define 
	$$\uppermu_n(\calS):=\limsup_{B_1,\dots,B_{n-1} \rightarrow \infty} 
		\limsup_{B_n \rightarrow \infty} 
		\frac{\#(\calS \cap \BOX)}{\# \BOX}.$$
This has the effect of considering only boxes
in which the $n^{\operatorname{th}}$ dimension is large relative to the others.
Define $\lowermu_n(\calS)$ similarly.
If $\uppermu_n(\calS)=\lowermu_n(\calS)$, define $\mu_n(\calS)$ as the
common value.
Also define
	$$\uppermu_\weak(\calS):=\max_{\sigma} 
		\limsup_{B_{\sigma(1)} \rightarrow\infty}
		\cdots
		\limsup_{B_{\sigma(n)} \rightarrow\infty}
		\frac{\#(\calS \cap \BOX)}{\# \BOX},$$
where $\sigma$ ranges over permutations of $\{1,2,\ldots,n\}$.
This definition in effect considers only boxes
whose dimensions can be ordered so that each is very large relative
to the previous dimensions.
Define $\lowermu_\weak(\calS)$ similarly,
and define $\mu_\weak(\calS)$ 
if $\uppermu_\weak(\calS)=\lowermu_\weak(\calS)$.

\section{Theorems}
\label{theorems}

Throughout this paper, $p$ represents a prime number.
In particular, in a sum or product indexed by $p$,
it is assumed that $p$ runs through only primes.
Similarly, $\pp$ represents a nonzero prime of $A$.

\begin{theorem}[Relatively prime values]
\label{relativelyprime}
Let $f,g \in A[x_1,\dots,x_n]$ be polynomials that are
relatively prime as elements of $K[x_1,\dots,x_n]$.
Let
	$$\calR_{f,g}:=\{\, a\in A^n : \gcd(f(a),g(a))=1 \,\}.$$
Then $\mu(\calR_{f,g})=\prod_{\pp} (1 - c_\pp/|\pp|^n)$,
where $\pp$ ranges over all nonzero primes of $A$,
and $c_\pp$ is the number of $x \in (A/\pp)^n$
satisfying $f(x)=g(x)=0$ in $A/\pp$.
\end{theorem}

The assumptions and conclusions for the squarefree value theorem
differ slightly in the $\Z$ and $\F_q[t]$ cases,
so we separate them into Theorem~\ref{Zsquarefree}
and Theorem~\ref{Fqtsquarefree}.

\begin{theorem}[Squarefree values over $\Z$]
\label{Zsquarefree}
Assume the $abc$ conjecture.
Let $f \in \Z[x_1,\dots,x_n]$ be a polynomial that is
squarefree as an element of $\Q[x_1,\dots,x_n]$,
and suppose that $x_n$ appears in $f$.
Let
	$$\calS_f:=\{\, a\in \Z^n : \text{$f(a)$ is squarefree} \,\}.$$
For each prime $p$,
let $c_p$ be the number of $x \in (\Z/p^2)^n$
satisfying $f(x)=0$ in $\Z/p^2$.
Then $\mu_n(\calS_f)=\prod_{p} (1 - c_p/p^{2n})$.
\end{theorem}

If the degree of $x_n$ in each irreducible factor of $f$ 
in Theorem~\ref{Zsquarefree} is $\le 3$,
then it is unnecessary to assume the $abc$ conjecture,
because the proof reduces to the case of one-variable polynomials
of degree $\le 3$, 
for which an unconditional result is known~\cite{hooley1967}.

The $n=1$ case of Theorem~\ref{Zsquarefree} differs slightly
from Theorem~1 in~\cite{granville1998}
in that the latter computes the density
of squarefree values of $f(x)/m$
where $m$ is a particular positive integer 
dividing all values of $f$.
Such results can be proved in the multivariable case
just as easily as Theorem~\ref{Zsquarefree}; 
the key to all such results is Lemma~\ref{squarefreeerror}.

\begin{cor}
\label{Zcorollary}
Let notation and assumptions be as in Theorem~\ref{Zsquarefree}
but without the restriction that $x_n$ appears in $f$.
Then $\mu_\weak(\calS_f)=\prod_{p} (1 - c_p/p^{2n})$.
\end{cor}

Corollary~\ref{Zcorollary} for $f(x_1,\dots,x_n)$
follows from Theorem~\ref{Zsquarefree}
applied to $f(x_{\sigma(1)},\dots,x_{\sigma(n)})$ for all 
permutations $\sigma$,
since in the definition of $\mu_\weak$ 
we may discard each $\limsup$ corresponding to 
a variable that does not appear.

\begin{theorem}[Squarefree values over $\Fqt$]
\label{Fqtsquarefree}
Let $A=\F_q[t]$.
Let $f \in A[x_1,\dots,x_n]$ be a polynomial that is
squarefree as an element of $K[x_1,\dots,x_n]$.
Let
	$$\calS_f:=\{\, a\in A^n : \text{$f(a)$ is squarefree} \,\}.$$
For each nonzero prime $\pp \subseteq A$,
let $c_\pp$ be the number of $x \in (A/\pp^2)^n$
satisfying $f(x)=0$ in $A/\pp^2$.
Then $\mu(\calS_f)=\prod_{\pp} (1 - c_\pp/|\pp|^{2n})$.
\end{theorem}

\begin{rem}
Note in particular that Theorem~\ref{Fqtsquarefree} proves a result
for $\mu$ instead of only for $\mu_n$.
\end{rem}

Suppose $\gcd(q,6)=1$.
One application of Theorem~\ref{Fqtsquarefree}
is to computing asymptotics for the weighted number $R_d$
of isomorphism classes of elliptic curves $(E,O)$ over $\F_q[t]$
with squarefree discriminant, as $d \rightarrow \infty$
for fixed $q$~\cite[\S4.22]{dejong2001}.
``Weighted'' means that each isomorphism class
receives the weight $1/\#\Aut(E,O)$ instead of $1$.
This number $R_d$ is closely connected
to the density of $(A,B) \in \F_q[t]^2$
such that the discriminant $\Delta=-16(4A^3+27B^2)$
of $y^2=x^3+Ax+B$ is squarefree,
except that one works with homogeneous polynomials 
$A \in H^0(\PP^1,\OO(4d))$
and $B \in H^0(\PP^1,\OO(6d))$,
so that the density has 
a factor corresponding to the point at infinity on $\PP^1$
in addition to the affine points.
A calculation shows that
the density of such $(A,B)$ having a double zero at a particular
closed point $\pp$ of $\PP^1=\PP^1_{\F_q}$ is $(2|\pp|^2-|\pp|)/|\pp|^4$,
where $|\pp|$ denotes the size of the residue field of $\pp$;
from this and our methods we obtain
	$$\lim_{d \rightarrow \infty} \frac{R_d}{q^{10d+1}} 
	= \frac{q}{q-1} \prod_{\pp \in \PP^1} 
		\left(1 - \frac{2|\pp|^2-|\pp|}{|\pp|^4}\right),$$
and the number $\gamma_q$ of~\cite[\S4.22]{dejong2001}
equals
	$$\gamma_q = \frac{q^3}{(q-1)^2(q+1)} \prod_{\pp \in \PP^1} 
		\left(1 - \frac{2|\pp|^2-|\pp|}{|\pp|^4}\right).$$

\begin{rem}
Since $-16(4A^3+27B^2)$ has degree only $2$ in $B$,
it is possible to obtain the result of the previous paragraph
by arguments simpler than those needed for the proof of 
Theorem~\ref{Fqtsquarefree} in general.
(This was pointed out to me by de Jong.)
\end{rem}

For a generalization of Theorem~\ref{Fqtsquarefree} 
to other rings of functions, 
see Section~\ref{ringsoffunctions}.
Analogues where ``squarefree'' is replaced 
by ``$k^{\operatorname{th}}$-power-free''
follow immediately from the same arguments
once one has Lemma~\ref{squarefreeerror}
(or the corresponding function field result).

\begin{theorem}[Values in $\Q^*/\Q^{*2}$]
\label{Qsquaredtheorem}
Let $f(x) \in \Z[x]$ be a nonzero polynomial.
Write $f(x) = c g(x)^2 h(x)$
where $c \in \Z$, $g(x) \in \Z[x]$,
and $h(x)$ is a squarefree polynomial in $\Z[x]$
whose coefficients have gcd~$1$.
If $\deg h>3$, assume the $abc$ conjecture.
Then the image of $\{f(1),f(2),\dots,f(B)\}$
in $\left( \Q^*/\Q^{*2} \right) \cup \{0\}$
has size $c_f B + o(B)$ for some constant $c_f \in [0,1]$.
If $\deg h=0$, then $c_f=0$.
If $\deg h\ge 2$, then $c_f=1$.
If $\deg h=1$, say $h(x)=ax+b$,
then 
	$$c_f = \frac{6}{\pi^2} \left( \sum_{r=0}^{|a|-1} \delta_r  \right)
	\prod_{p \mid a} (1-p^{-2})^{-1} \quad \in \frac{1}{\pi^2} \Q,$$
where $\delta_r :=1/m^2$ if $m$ is the smallest positive integer
satisfying $m^2 r \equiv b \pmod{a}$,
or $\delta_r:=0$ if no such $m$ exists.
\end{theorem}

Assuming the $abc$ conjecture, 
Granville~\cite[Corollary~2]{granville1998} proved
that the size of the image in Theorem~\ref{Qsquaredtheorem}
was {\em at least} some positive constant
times $B$ (when $f(x)$ has no repeated roots),
and guessed that the size should be asymptotic
to a constant times $B$, 
as our Theorem~\ref{Qsquaredtheorem} shows.
An essentially equivalent version of Theorem~\ref{Qsquaredtheorem} 
has been independently proved 
by P.~Cutter, A.~Granville, and T.~Tucker
(Theorems 1A, 1B, and 1C of~\cite{cuttergranvilletucker}),
using a similar proof.
They also prove a few related results not considered here.

It is natural to ask, as Granville has also done,
what the multivariable analogue of Theorem~\ref{Qsquaredtheorem} 
should be.
Here we formulate a precise question along these lines:

\begin{question}
Suppose $f \in \Z[x_1,\dots,x_n]$ is nonconstant and squarefree
as an element of $\Q[x_1,\dots,x_n]$.
For $B \ge 1$, let $S_B = f(\{1,2,\dots,B\}^n) \subset \Z$,
and let $T_B$ be the image of $S_B$ 
in $\left( \Q^*/\Q^{*2} \right) \cup \{0\}$.
Does $\#T_B/\#S_B$ tend to a positive limit as $B \rightarrow \infty$?
\end{question}

We do not have enough evidence to conjecture an answer.
But even if the answer is yes, it is not clear that we would understand
the asymptotic size of $T_B$, 
because even the problem of estimating $\#S_B$ seems very difficult.

\section{Zero values}
\label{zerosection}

The following lemma is well-known.
We include a proof mainly because it is a toy version
of some of the reductions used later on.

\begin{lemma}
\label{zeros}
Let $f \in A[x_1,\dots,x_n]$ be a nonzero polynomial.
Let $\calZ = \{\, a\in A^n : f(a)=0 \,\}$.
Then $\mu(\calZ)=0$.
\end{lemma}

\begin{proof}
We use induction on $n$.
The base case $n=0$ is trivial, so suppose $n \ge 1$.
Let $f_1 \in A[x_1,\ldots,x_{n-1}]$
be the leading coefficient of $f$
when $f$ is viewed as a polynomial in $x_n$.
Let $\delta$ be the $x_n$-degree of $f$.
Now $\calZ \subseteq \calZ_1 \cup \calZ_2$ where
\begin{align*}
	\calZ_1 &:= \{\, a\in A^n : f_1(a)=0 \,\}, \\
	\calZ_2 &:= \{\, a \in A^n: f_1(a) \not=0 \text{ and } f(a)=0 \,\}.
\end{align*}
By the inductive hypothesis, $\mu(\calZ_1)=0$.
For each $(a_1,\dots,a_{n-1}) \in A^{n-1}$,
there are at most $\delta$ values $a_n \in A$
for which $(a_1,\dots,a_{n-1},a_n) \in \calZ_2$.
Thus $\mu(\calZ_2)=0$, by definition of $\mu$.
Hence $\mu(\calZ)=0$, as desired.
\end{proof}

\section{Relatively prime values}
\label{relativelyprimesection}

The bulk of the work in proving Theorem~\ref{relativelyprime}
is in the following.

\begin{lemma}
\label{relativelyprimeerror}
Let $f,g \in A[x_1,\dots,x_n]$ be polynomials that are
relatively prime as elements of $K[x_1,\dots,x_n]$.
Let
	$$\calQ_{f,g,M}:=\{\, a\in A^n : \exists \pp 
	\text{ such that $|\pp| \ge M$ and $\pp \mid f(a),g(a)$} \,\}.$$
Then $\lim_{M \rightarrow \infty} \uppermu(\calQ_{f,g,M})=0$.
\end{lemma}

\begin{proof}
Since we are interested only in $\pp$ with $|\pp|$ large,
we may divide $f$ and $g$ by any factors in $A$ that they have,
in order to assume that $f$ and $g$ are 
relatively prime as elements of $A[x_1,\dots,x_n]$.

The proof will be by induction on $n$.
The case $n=0$ is trivial, so assume $n \ge 1$.
We need to bound the size of $Q := \calQ_{f,g,M} \cap \BOX$,
whenever the ``dimensions'' $B_i$ of $\BOX$ are sufficiently large.
Without loss of generality, $M \le B_1 \le B_2 \le \cdots \le B_n$.
Set $B_0=M$ and $B_{n+1}=\infty$.
Let $f_1,g_1 \in A[x_1,\ldots,x_{n-1}]$
be the leading coefficients of $f$ and $g$
when $f$ and $g$ are viewed as polynomials in $x_n$.

\medskip
\noindent{\em Case 1.} One of the polynomials, say $g$, is
a polynomial in $x_1,\dots,x_{n-1}$ only.

In this case, we use an inner induction on $\delta$,
where $\delta$ is the $x_n$-degree of $f$.
The base case $\delta=0$ is handled by the outer inductive hypothesis,
so from now on assume $\delta>0$.
We may reduce to the case that $f$ and $g$ are irreducible.
If $g \mid f_1$, then we can subtract a multiple of $g$ from $f$
to lower its $x_n$-degree $\delta$, without changing
$\calQ_{f,g,M}$ or the relative primality of $f$ and $g$,
so the result follows from the inner inductive hypothesis.
Hence we may assume $g \nmid f_1$.
Since $g$ is irreducible, $f_1$ and $g$ are relatively prime
in $A[x_1,\dots,x_{n-1}]$.

Now $Q=\bigcup_{s=0}^n Q_s$,
where 
	$$Q_s :=\{\, a \in \BOX : \exists \pp 
	\text{ such that $B_s \le |\pp| < B_{s+1}$ 
		and $\pp \mid f(a),g(a)$} \,\},$$
so it suffices to show that given $0 \le s \le n$,
the ratio $\#Q_s/\#\BOX$ can be made arbitrarily small
by choosing the $B_i$ sufficiently large.

Suppose we fix $s$ with $0 \le s <n$.
(We will bound $Q_n$ later.)
Let $X$ be the subscheme of $\Aff^n_A$ defined by $f=g=0$.
Since $f$ and $g$ are relatively prime, $X$ has codimension at least $2$
in $\Aff^n_A$.
Let $\pi: \Aff^n_A \rightarrow \Aff^s_A$ be the projection onto the
first $s$ coordinates.
Let $Y_i$ be the (constructible) set 
of $y \in \Aff^s_A$ such that the fiber $X_y:= X \cap \pi^{-1}(y)$
has codimension $i$ in $\pi^{-1}(y) \isom \Aff^{n-s}_{\kappa(y)}$.
(Here $\kappa(y)$ denotes the residue field of $y$.)
Since $X$ has codimension at least $2$ in $\Aff^n_A$,
it follows from Theorem~15.1(i) of~\cite{matsumura1989}
that the subset $Y_i$ has codimension at least $2-i$ in $\Aff^s_A$.
In particular, we can choose a nonzero $h \in A[x_1,\dots,x_s]$
vanishing on $Y_1$.
Also we can find
relatively prime $j_1,j_2 \in A[x_1,\dots,x_s]$ vanishing on $Y_0$ 
as follows: choose any nonzero $j_1$ vanishing on $Y_0$;
if $I(Y_0)$ is not contained in the union of the minimal primes
over $(j_1)$, then any $j_2 \in I(Y_0)$ outside those primes will be
relatively prime to $j_1$; if $I(Y_0)$ is contained in that union,
then Proposition~1.11(i) of~\cite{atiyah-macdonald} implies that
$I(Y_0)$ is contained in some minimal prime over $(j_1)$,
but such a prime has codimension~$1$,
contradicting the fact that $Y_0$ has codimension at least $2$.
Define $Y_{\ge 2} := \bigcup_{i \ge 2} Y_i$.

Given $a=(a_1,\dots,a_n) \in A^n$ and a nonzero prime $\pp$ of $A$,
let $a_\pp = (a_1,\dots,a_n)_\pp$ denote the closed point in $\Aff^n_{A/\pp}$
whose coordinates are $a_1,\dots,a_n$.
Thus 
	$$Q_s = \{\, a \in \BOX : \exists \pp 
	\text{ such that $B_s \le |\pp| < B_{s+1}$ 
		and $a_\pp \in X$} \,\}.$$
Let $Z:=\{\, a \in \BOX : h(a_1,\dots,a_s)=0  \,\}$.
Define
	$$R_{\ge 2} := \{\, a \in \BOX : 
		\text{ $\exists \pp$ such that 
		$B_s \le |\pp| < B_{s+1}$, $a_\pp \in X$, 
		and $(a_1,\dots,a_s)_\pp \in Y_{\ge 2}$} \,\},$$
and define $R_1$ and $R_0$ similarly, using $Y_1$ and $Y_0$, respectively,
in place of $Y_{\ge 2}$.

Then $Q_s \subseteq Z \cup R_{\ge 2} \cup (R_1-Z) \cup R_0$.
By Lemma~\ref{zeros}, $\#Z/\#\BOX$ can be made arbitrarily small
by choosing the $B_i$ sufficiently large.

Next consider $R_{\ge 2}$.
It suffices to show that for $(a_1,\dots,a_s) \in \BOX(B_1,\dots,B_s)$,
the fraction of $(a_{s+1},\dots,a_n)$ in $\BOX(B_{s+1},\dots,B_n)$
for which there exists a prime $\pp$
with $B_s \le |\pp| < B_{s+1}$, $a_\pp \in X$, 
and $(a_1,\dots,a_s)_\pp \in Y_{\ge 2}$
is small when $B_s$ is large.
Fix $(a_1,\dots,a_s) \in \BOX(B_1,\dots,B_s)$.
If $\pp$ is a prime with $B_s \le |\pp| < B_{s+1}$
and $y:=(a_1,\dots,a_s)_\pp$ lies in $Y_{\ge 2}$,
then $X_y$ has codimension at least $2$ in $\Aff^{n-s}_{A/\pp}$,
so 
	$$\# X_y(A/\pp) = O(|\pp|^{n-s-2}) = O((\#(A/\pp)^{n-s}) / |\pp|^2).$$
Moreover, the implied constant can be made uniform in $y$,
since the $X_y$ are fibers in an algebraic family.
Since $|\pp|<B_{s+1}$,
the reductions modulo $\pp$
of the $(a_{s+1},\dots,a_n) \in \BOX(B_{s+1},\dots,B_n)$
are almost uniformly distributed in $(A/\pp)^{n-s}$:
to be precise,
each residue class in $(A/\pp)^{n-s}$ 
is represented by a fraction at most $O(\#(A/\pp)^{-(n-s)})$ 
of these $(a_{s+1},\dots,a_n)$,
where the implied constant depends only on $n$.
Hence the fraction of $(a_{s+1},\dots,a_n) \in \BOX(B_{s+1},\dots,B_n)$
satisfying $(a_1,\dots,a_n)_\pp \in X_y$
is $O(1/|\pp|^2)$,
and summing over all $\pp$ with $B_s \le |\pp| < B_{s+1}$
still yields a fraction that can be made arbitrarily small
by taking $B_s$ large, since $\sum_\pp 1/|\pp|^2$ converges.

We now adapt the previous paragraph to bound $\#(R_1-Z)$.
Suppose $(a_1,\dots,a_s) \in \BOX(B_1,\dots,B_s)$
and $h(a_1,\dots,a_s)\not=0$.
Let $\eta$ be the total degree of $h$.
Then $|h(a_1,\dots,a_s)| = O(B_s^{\eta})$,
where the constant implied by the $O$ depends only on $h$,
not on the $a_i$ or $B_i$.
Thus, provided that $B_s$ is large,
$h(a_1,\dots,a_s)$ can be divisible by at most $\eta$
primes $\pp$ satisfying $B_s \le |\pp| < B_{s+1}$.
Hence $(a_1,\dots,a_s)_\pp \in Y_1$ for at most $\eta$
primes $\pp$ satisfying $B_s \le |\pp| < B_{s+1}$.
By definition of $Y_1$, if $y=(a_1,\dots,a_s)_\pp$ for such $\pp$,
then $X_y$ has codimension at least $1$ in $\Aff^{n-s}_{A/\pp}$,
so $\# X_y(A/\pp) = O((\#(A/\pp)^{n-s}) / |\pp|)$,
where the implied constant is independent of $y$.
The reductions modulo $\pp$
of the $(a_{s+1},\dots,a_n) \in \BOX(B_{s+1},\dots,B_n)$
are again almost uniformly distributed in $(A/\pp)^{n-s}$.
Hence the fraction of $(a_{s+1},\dots,a_n)$ in $\BOX(B_{s+1},\dots,B_n)$
whose reduction modulo $\pp$ lies in $X_y$
is $O(1/|\pp|)$.
Summing over at most $\eta$ possible primes $\pp$ with $|\pp| \ge B_s$
still yields a fraction that can be made arbitrarily small
by taking $B_s$ large.

Finally we consider $R_0$.
Since $s<n$, the outer inductive hypothesis applied to $j_1$ and $j_2$
implies that $\#R_0/\#\BOX$ can be made arbitrarily small
by taking the $B_i$ large.

To finish Case~1,
we need to bound $Q_n$.
We have $Q_n \subseteq S_0 \cup S \cup S'$, where 
\begin{align*}
	S_0 &:=\{\, a \in \BOX : g(a_1,\dots,a_{n-1})=0  \,\}, \\
	S &:= \{\, a \in \BOX : \exists \pp \text{ such that 
			$|\pp| \ge B_n$ and $\pp \mid f_1(a),g(a)$} \}\, \\
	S' &:= \{\, a \in \BOX : g(a_1,\dots,a_{n-1}) \not= 0 
		\text{ and $\exists \pp$ such that 
		$|\pp| \ge B_n$,
		$\pp \mid f(a),g(a)$ 
		and $\pp \nmid f_1(a)$} \,\}.
\end{align*}
Lemma~\ref{zeros} bounds $\#S_0/\#\BOX$.
The outer inductive hypothesis applied to $f_1$ and $g$
bounds $\#S/\#\BOX$.

It remains to bound $\#S'$.
For $(a_1,\dots,a_{n-1}) \in \BOX(B_1,\dots,B_{n-1})$
such that $g(a_1,\dots,a_{n-1}) \not=0$,
we will show that the fraction of $a_n \in \BOX(B_n)$
such that there exists $\pp$ with
$|\pp| \ge B_n$, $\pp \mid f(a),g(a)$ and $\pp \nmid f_1(a)$
is small.
We use a method similar to that used to bound $R_1$.
Let $\gamma$ denote the total degree of $g$.
If $B_n$ is sufficiently large (depending only on $g$),
then given $(a_1,\dots,a_{n-1})$,
there are at most $\gamma$ primes $\pp$ dividing $g(a)$
with $|\pp|\ge B_n$.
For each such $\pp$, if moreover $\pp \nmid f_1(a)$,
then the polynomial $f(a_1,\dots,a_{n-1},x_n) \bmod \pp$
in $(A/\pp)[x_n]$ is of degree $\delta$,
and has at most $\delta$ roots in $A/\pp$.
For each such root, there are at most $O(1)$ elements of $\BOX(B_n)$
reducing to it modulo $\pp$, since $|\pp| \ge B_n$.
Thus given $(a_1,\dots,a_{n-1})$,
there are at most $\gamma \delta \cdot O(1) = O(1)$ 
values of $a_n \in \BOX(B_n)$ for which $(a_1,\dots,a_n) \in S'$.
Thus $\#S'/\#\BOX$ can be made small by choosing $B_n$ large.

\medskip
\noindent{\em Case 2.} The $x_n$-degree of $f$ and $g$ are both positive.

Let $R \in A[x_1,\ldots,x_{n-1}]$
be the resultant of $f$ and $g$ with respect to $x_n$.
Since $f$ and $g$ are relatively prime, $R$ is nonzero.
Since $f_1$, $g_1$, and $R$ are all nonzero and do not involve $x_n$,
none of them are multiples of $f$ or $g$.
Since $f$ and $g$ are irreducible,
each of $f_1$, $g_1$, and $R$ must be relatively prime to each of $f$ and $g$.
Moreover, if $\pp$ is a prime dividing $f(a)$ and $g(a)$,
and if the leading coefficients $f_1(a)$ and $g_1(a)$ are nonzero modulo $\pp$,
then by a well known property of the resultant, $\pp \mid R(a)$.
Hence
\begin{align*}
	\{\, a \in \BOX: \pp \mid f(a), g(a) \,\}
	&\subseteq \{\, a \in \BOX: \pp \mid f_1(a), g(a) \,\} \\
	&\qquad	\cup \{\, a \in \BOX: \pp \mid f(a), g_1(a) \,\} \\
	&\qquad	\cup \{\, a \in \BOX: \pp \mid f(a), R(a) \,\}.
\end{align*}
Taking the union over all $\pp$ with $|\pp| \ge M$ 
and applying Case~1 to $f_1,g$, to $f,g_1$, and to $f,R$
completes the proof.
\end{proof}

\begin{proof}[Proof of Theorem~\ref{relativelyprime}]
Let $P_M$ denote the set of nonzero primes $\pp$ of $A$ such that $|\pp| < M$.
Approximate $\calR_{f,g}$ by
	$$\calR_{f,g,M} := \{\, a\in A^n : 
		\text{ $f(a)$ and $g(a)$ are not both divisible 
		by any prime $\pp \in P_M$} \,\}.$$
Define the ideal $I$ as the product of all $\pp$ in $P_M$.
Then $\calR_{f,g,M}$ is a union of cosets of 
the subgroup $I^n \subset A^n$.
(Here $I^n$ is the cartesian product.)
Hence $\mu(\calR_{f,g,M})$ is the fraction of residue classes in $(A/I)^n$
in which for all $\pp \in P_M$,
at least one of $f(a)$ and $g(a)$ is nonzero modulo $\pp$.
Applying the Chinese Remainder Theorem shows that
$\mu(\calR_{f,g,M}) = \prod_{\pp \in P_M} (1 - c_{\pp} / |\pp|^n)$.
By Lemma~\ref{relativelyprimeerror},
	$$\mu(\calR_{f,g}) = 
	\lim_{M \rightarrow \infty} \mu(\calR_{f,g,M}) 
	= \prod_{\pp} (1 - c_{\pp} / |\pp|^n).$$
Since $f$ and $g$ are relatively prime as elements of $K[x_1,\dots,x_n]$,
there exists a nonzero $u \in A$ such that
$f=g=0$ defines a subscheme of $\Aff^n_{A[1/u]}$ of codimension
at least $2$.
Thus $c_{\pp} = O(|\pp|^{n-2})$ as $|\pp| \rightarrow \infty$, 
and the product converges.
\end{proof}

\section{Squarefree values of polynomials over $\Z$}
\label{Zsection}

If $f \in A[x_1,\dots,x_n]$, and $M \ge 1$,
define
	$$\calT_{f,M}:=\{\, a\in A^n : \text{$\exists \pp$ with $|\pp| \ge M$
	such that $\pp^2 \mid f(a)$} \,\}.$$
For the rest of this section, we take $A=\Z$.
The following is a variant of Theorem~1 of~\cite{granville1998},
and has the same proof.

\begin{lemma}
\label{squarefreeerror1}
Assume the $abc$ conjecture.
Suppose that $f \in \Z[x]$
is squarefree as a polynomial in $\Q[x]$.
For each prime $p \ge M$, let
$c_p$ be the number of $x \in \Z/p^2$ satisfying $f(x) = 0$ in $\Z/p^2$.
Then $1 - \mu(\calT_{f,M})= \prod_{p \ge M} (1-c_p/p^2)$.
\end{lemma}

We are now ready to prove 
the analogue of Lemma~\ref{relativelyprimeerror}
for squarefree values of multivariable polynomials over $\Z$:

\begin{lemma}
\label{squarefreeerror}
Assume the $abc$ conjecture.
Suppose that $f \in \Z[x_1,\dots,x_n]$ 
is squarefree as a polynomial in $\Q[x_1,\dots,x_n]$,
and suppose that $x_n$ appears in $f(x)$.
Then $\lim_{M \rightarrow \infty} \uppermu_n(\calT_{f,M})=0$.
\end{lemma}

\begin{proof}
Factors of $f$ lying in $\Z$ are irrelevant as $M \rightarrow \infty$,
so we may assume that $f$ is squarefree as 
a polynomial in $\Z[x_1,\dots,x_n]$.
If $f$ factors as a product of two relatively prime polynomials
$g$ and $h$, then the result for $f$ follows from the result
for $g$ and $h$ together with Lemma~\ref{relativelyprimeerror}
applied to $g,h$.
Hence we may reduce to the case where $f$ is irreducible
in $\Z[x_1,\dots,x_n]$.

Let $\Delta \in \Z[x_1,\dots,x_{n-1}]$ and $\delta \ge 1$ 
be the discriminant and degree, respectively,
of $f$ considered as a polynomial in $x_n$.
Given $B_1,\dots,B_n$, let $Q := \calT_{f,M} \cap \BOX$.
We need to show that if the $B_i$ are sufficiently large,
and $B_n$ is sufficiently large relative to the other $B_i$,
then $\#Q/\#\BOX$ is small.

The fraction of $(a_1,\dots,a_{n-1})$
in $\BOX_{n-1}:=\BOX(B_1,\dots,B_{n-1})$
at which $\Delta$ vanishes is negligible, by Lemma~\ref{zeros}.
Since $f$ is irreducible,
when $f$ is viewed as a polynomial in $x_n$,
its coefficients (in $\Z[x_1,\dots,x_{n-1}]$)
are relatively prime (not necessarily pairwise).
In particular, the common zero locus of these coefficients
has codimension at least $2$ in $\Aff^{n-1}_\Q$, 
hence is contained in the subvariety defined by $\tilde{f}=\tilde{g}=0$
for two relatively prime elements 
$\tilde{f}, \tilde{g} \in \Z[x_1,\dots,x_{n-1}]$.
Thus Lemma~\ref{relativelyprimeerror} implies
that the fraction of $(a_1,\dots,a_{n-1}) \in \BOX_{n-1}$
such that there exists a prime $p \ge M$
such that the image of $f(a_1,\dots,a_{n-1},x_n)$ in $\F_p[x_n]$
is zero is negligible, when $M$ is large.

It remains to bound $\#(Q \cap (Q' \times [1,B_n]))/\#\BOX$,
where $Q'$ is the set of $(a_1,\dots,a_{n-1}) \in \BOX_{n-1}$
such that
\begin{itemize}
\item $\Delta(a_1,\dots,a_{n-1}) \not= 0$, and
\item there is no prime $p \ge M$ such that 
	the image of $f(a_1,\dots,a_{n-1},x_n)$ in $\F_p[x_n]$ is zero.
\end{itemize}
By the first condition, Lemma~\ref{squarefreeerror1}
applies to $f(a_1,\dots,a_{n-1},x_n) \in \Z[x_n]$
for each $(a_1,\dots,a_{n-1}) \in Q'$.
Letting $B_n$ tend to infinity while $B_1,\dots,B_{n-1}$ are fixed,
we find that it suffices to bound
\begin{equation}
\label{densityaverage}
	\frac{1}{\#\BOX} \sum_{(a_1,\dots,a_{n-1}) \in Q'} 
		B_n \left( 1 - \prod_{p \ge M} 
		\left( 1 - \frac{c_p(a_1,\dots,a_{n-1})}{p^2} \right) \right),
\end{equation}
where $c_p(a_1,\dots,a_{n-1})$ is the number of $x_n \in \Z/p^2$
such that $f(a_1,\dots,a_{n-1},x_n)=0$ in $\Z/p^2$.
The inequality $1-\alpha \beta \le (1-\alpha) + (1-\beta)$
holds for $\alpha,\beta \in [0,1]$;
applying this with
	$$\alpha := 
	\prod_{\twoline{p \ge M}{p \nmid \Delta(a_1,\dots,a_{n-1})}}
	   \left( 1 - \frac{c_p(a_1,\dots,a_{n-1})}{p^2} \right), \qquad
	\beta := \prod_{\twoline{p \ge M}{p \mid \Delta(a_1,\dots,a_{n-1})}}
	   \left( 1 - \frac{c_p(a_1,\dots,a_{n-1})}{p^2} \right)$$
and using
	$$1-\alpha \le 
	\sum_{\twoline{p \ge M}{p \nmid \Delta(a_1,\dots,a_{n-1})}} 
		\frac{c_p(a_1,\dots,a_{n-1})}{p^2}$$
bounds~(\ref{densityaverage}) by $s_1+s_2$ where
\begin{align*}
	s_1 &:= \frac{1}{\#\BOX} \sum_{(a_1,\dots,a_{n-1}) \in Q'} 
		B_n 
		\sum_{\twoline{p \ge M}{p \nmid \Delta(a_1,\dots,a_{n-1})}} 
		\frac{c_p(a_1,\dots,a_{n-1})}{p^2}, \\
	s_2 &:= \frac{1}{\#\BOX} \sum_{(a_1,\dots,a_{n-1}) \in Q'} 
		B_n \left( 1 - 
		\prod_{\twoline{p \ge M}{p \mid \Delta(a_1,\dots,a_{n-1})}}
		\left( 1 - \frac{c_p(a_1,\dots,a_{n-1})}{p^2} \right) \right).
\end{align*}
When $(a_1,\dots,a_{n-1}) \in Q'$ and $p \nmid \Delta(a_1,\dots,a_{n-1})$,
Hensel's Lemma implies $c_p(a_1,\dots,a_{n-1}) \le \delta$,
while $B_n \#Q' \le \#\BOX$, so $s_1 \le \sum_{p \ge M} \delta/p^2$,
which is negligible as $M \rightarrow \infty$.
When $(a_1,\dots,a_{n-1}) \in Q'$ and $p \mid \Delta(a_1,\dots,a_{n-1})$,
the image of $f(a_1,\dots,a_{n-1},x_n)$ in $\F_p[x_n]$
has at most $\delta$ zeros in $\F_p$;\,
so $c_p(a_1,\dots,a_{n-1}) \le \delta p$,
and
\begin{equation}
\label{s2inequality}
	s_2 \le \frac{1}{\#\BOX} \sum_{(a_1,\dots,a_{n-1}) \in Q'} 
		B_n \left( 1 - 
		\prod_{\twoline{p \ge M}{p \mid \Delta(a_1,\dots,a_{n-1})}}
		\left( 1 - \frac{\delta}{p} \right) \right).
\end{equation}
Let $\Phi(x_n)=\prod_{j=1}^\delta (x_n-j)$.
We may assume $M \ge \delta$; then
	$$1 - \prod_{\twoline{p \ge M}{p \mid \Delta(a_1,\dots,a_{n-1})}}
		\left( 1 - \frac{\delta}{p} \right)$$
equals the density of 
	$$\{\, x_n \in \Z : \text{ $\exists p \ge M$ such that 
			$p \mid \Delta(a_1,\dots,a_{n-1}),
			\Phi(x_n)$} \,\}.$$
Hence, as $B_n \rightarrow \infty$ for fixed $M,B_1,\dots,B_{n-1}$,
the right hand side of~(\ref{s2inequality}) has the same limit as
	$$\frac{\# ( (Q' \times [1,B_n]) \cap 
		\calQ_{\Delta(x_1,\dots,x_{n-1}),\Phi(x_n),M} )}
	{\#\BOX}.$$
The latter is negligible, by Lemma~\ref{relativelyprimeerror}.
\end{proof}

\begin{rem}
It seems difficult to improve Lemma~\ref{squarefreeerror} to obtain a result
for the more natural definition of density, 
$\uppermu$ instead of $\uppermu_n$.
This would require a version of Granville's one-variable result that
is uniform in the coefficients of the polynomial.
Granville's proof uses Belyi functions, however, whose degrees
vary wildly with the coefficients.
\end{rem}

\begin{proof}[Proof of Theorem~\ref{Zsquarefree}]
Approximate $\calS_f$ by
	$$\calS_{f,M} := \{\, a\in \Z^n : 
		\text{ $f(a)$ is not divisible by $p^2$
		for any prime $p<M$} \,\}.$$
Then $\calS_{f,M}$ is a union of cosets of $(I\Z)^n$
where $I=\prod_{p<M} p^2$.
The Chinese Remainder Theorem implies
$\mu_n(\calS_{f,M}) = \prod_{\pp \in P_M} (1 - c_p / p^{2n})$.
By Lemma~\ref{squarefreeerror},
	$$\mu_n(\calS_f) = 
	\lim_{M \rightarrow \infty} \mu_n(\calS_{f,M}) 
	= \prod_p (1 - c_p / p^{2n}).$$

Finally, we show that the product converges (instead of diverging to $0$)
by showing that $c_p = O(p^{2n-2})$.
(The value of the product could still be zero if some factor were zero.)
Let $X$ be the subscheme of $\Aff^n_\Z$ defined by $f=0$.
Since the field $\Q$ is perfect,
the nonsmooth locus of $X \times \Q \rightarrow \Spec \Q$ 
has codimension at least $2$ in $\Aff^n_\Q$.
It follows that for sufficiently large $p$, 
the nonsmooth locus $Y_p$ of $X \times \F_p \rightarrow \Spec \F_p$ 
has codimension at least $2$ in $\Aff^n_{\F_p}$,
so $\#Y_p(\F_p) = O(p^{n-2})$.
Each point in $Y_p(\F_p)$ can be lifted to an $n$-tuple in $(\Z/p^2)^n$
in $p^n$ ways.
On the other hand, $\#(X_p-Y_p)(\F_p) = O(p^{n-1})$
but the nonvanishing of some derivative modulo $p$
at a point in $(X_p-Y_p)(\F_p)$ implies that such a point
lifts to at most $p^{n-1}$ solutions to $f(x)=0$ in $(\Z/p^2)^n$.
Thus $c_p = p^n O(p^{n-2}) + p^{n-1} O(p^{n-1}) = O(p^{2n-2})$.
\end{proof}

\section{Squarefree values of polynomials over $\F_q[t]$}
\label{Fqtsection}

Throughout this section, $A=\F_q[t]$ where $q=p^e$.
Our goal is to prove Theorem~\ref{Fqtsquarefree}.
We begin by stating the analogue of Lemma~\ref{squarefreeerror}.

\begin{lemma}
\label{Fqtsquarefreeerror}
Suppose that $f \in A[x_1,\dots,x_n]$ 
is squarefree as a polynomial in $K[x_1,\dots,x_n]$.
Then $\lim_{M \rightarrow \infty} \uppermu(\calT_{f,M})=0$.
\end{lemma}

Before beginning the proof of Lemma~\ref{Fqtsquarefreeerror},
we state and prove two results that will be needed in its proof.

\begin{lemma}
\label{restrictionofscalars}
If $f \in K[x_1,\dots,x_n]$ is squarefree, then
	$$F:=f\left( y_0^p + t y_1^p + \cdots + t^{p-1} y_{p-1}^p , 
			x_2,x_3,\ldots,x_n \right) 
	\in K[y_0,\dots,y_{p-1},x_2,\dots,x_n]$$
is squarefree.
\end{lemma}

\begin{proof}
We work in $B:=K^{1/p}[y_0,\dots,y_{p-1},x_2,\dots,x_n]$,
where $K^{1/p} = \F_q(t^{1/p})$.
Define $u:=y_0 + t^{1/p} y_1 + \cdots + t^{(p-1)/p} y_{p-1}$.

We first show that $K^{1/p}[u] \cap K[y_0,\dots,y_{p-1}] = K[u^p]$.
Suppose $g = \sum \alpha_i u^i \in K^{1/p}[u] \cap K[y_0,\dots,y_{p-1}]$.
The coefficient of $y_0^i$ in $g$ is $\alpha_i$,
so $\alpha_i \in K$ for all $i$.
The coefficient of $y_0^{i-1} y_1$ in $g$ is $i t^{1/p} \alpha_i$;
this too is in $K$, so for all $i$, either $p \mid i$ or $\alpha_i=0$.
Thus $g \in K[u^p]$, as desired.

It follows that
\begin{equation}
\label{weirdintersection}
	K^{1/p}[u,x_2,\ldots,x_n] \cap K[y_0,\dots,y_{p-1},x_2,\dots,x_n] 
	= K[u^p,x_2,\dots,x_n].
\end{equation}
Suppose that $G^2 \mid F$ for some
$G \in K[y_0,\dots,y_{p-1},x_2,\dots,x_n] - K$.
Then $G^2 \mid F$ in $B$.
If we view $B$ as a polynomial ring over $K^{1/p}$
in algebraically independent
indeterminates $u,y_1,\dots,y_{p-1},x_2,\dots,x_n$ (by eliminating $y_0$),
then $F = f(u^p,x_2,\ldots,x_n)$ does not involve $y_1,\dots,y_{p-1}$,
so $G \in K^{1/p}[u,x_2,\ldots,x_n]$ too.
By~(\ref{weirdintersection}),
$G \in K[u^p,x_2,\dots,x_n]$.
If we write $G = g(u^p,x_2,\dots,x_n)$ where $g \in K[x_1,\dots,x_n]$,
then $g^2 \mid f$ and $g \not\in K$, contradicting the assumption
that $f$ is squarefree in $K[x_1,\dots,x_n]$.
\end{proof}

\begin{lemma}
\label{derivative}
Suppose $f \in K[x_1^p,\dots,x_n^p]$ is squarefree as
an element of $K[x_1,\dots,x_n]$.
Then $f$ and $\del f/\del t$ are relatively prime as elements
of $K[x_1,\dots,x_n]$.
\end{lemma}

\begin{proof}
The gcd $g$ of $f$ and $\del f/\del t$ in $K[x_1^p,\dots,x_n^p]$
equals their gcd in $K[x_1,\dots,x_n]$.
We may multiply $g$ by an element of $K^*$ to assume
that some coefficient of $g$ equals $1$.
Write $f=gh$.
Since $f$ is squarefree, $g$ and $h$ are relatively prime 
in $K[x_1,\dots,x_n]$.
But $g$ divides 
$\frac{\del f}{\del t} = g \frac{\del h}{\del t} + h \frac{\del g}{\del t}$,
so $g$ divides $\del g/\del t$.
The total degree of $\del g/\del t$ is less than or equal to that
of $g$, so $\del g/\del t = c g$ for some $c \in K$.
Then each coefficient $\gamma$ of $g$ satisfies
$\del \gamma/\del t = c \gamma$.
One of these coefficients is $1$,
so $c=0$.
Thus each coefficient $\gamma$ is in $K^p$,
so $g=G^p$ for some $G \in K[x_1,\dots,x_n]$.
But $g \mid f$, and $f$ is squarefree,
so $g \in K$.
Hence $f$ and $\del f/\del t$ are relatively prime.
\end{proof}

\begin{proof}[Proof of Lemma~\ref{Fqtsquarefreeerror}]
We may assume that $f$ is squarefree as an element of $A[x_1,\dots,x_n]$.
Define
	$$F:=f \left( 
		\sum_{j=0}^{p-1} t^j y_{1j}^p,
		\ldots,
		\sum_{j=0}^{p-1} t^j y_{nj}^p \right) 
	\in K[\ldots,y_{ij},\ldots]$$
By $n$ applications of Lemma~\ref{restrictionofscalars},
$F$ is squarefree.
Let $B_{ij} = (B_i/|t^j|)^{1/p}$.
As each $y_{ij}$ ranges 
through elements of $A$ satisfying $|y_{ij}| \le B_{ij}$,
the $n$-tuples 
	$$\left( 
		\sum_{j=0}^{p-1} t^j y_{1j}^p,
		\ldots,
		\sum_{j=0}^{p-1} t^j y_{nj}^p \right)$$
exhaust the elements of $\BOX$, with each element appearing once.
Hence Lemma~\ref{Fqtsquarefreeerror} for $f$
follows from Lemma~\ref{Fqtsquarefreeerror} for $F$.

Evaluation of $F$ at an element $a \in A^{np}$
commutes with formal application of $\del/\del t$,
since $F \in K[\ldots,y_{ij}^p,\ldots]$.
Thus, for a prime $\pp$ of $A$, we have $\pp^2 \mid F(a)$
if and only if $\pp$ divides $F(a)$ and $(\del F/\del t)(a)$.
Hence Lemma~\ref{Fqtsquarefreeerror} for $F$
follows from Lemma~\ref{relativelyprimeerror} for $F$ and $\del F/\del t$,
which are relatively prime by Lemma~\ref{derivative}.
\end{proof}

\begin{proof}[Proof of Theorem~\ref{Fqtsquarefree}]
We mimic the proof of Theorem~\ref{Zsquarefree} at the end of 
Section~\ref{Zsection}, 
using Lemma~\ref{Fqtsquarefreeerror} 
in place of Lemma~\ref{squarefreeerror}.
But the last paragraph of that proof, proving $c_p=O(p^{2n-2})$
to obtain convergence of the infinite product,
does not carry over, since it uses the fact that $\Q$ is perfect
to obtain generic smoothness, and $\F_q(t)$ is not perfect.

Therefore we prove $c_\pp = O(|\pp|^{2n-2})$ by a different method.
The fraction $c_\pp/|\pp|^{2n}$ is unchanged when we replace 
the $c_\pp$ and $n$ for the original $f$ by 
the corresponding values for the $F$ defined in 
the beginning of the proof of Lemma~\ref{Fqtsquarefreeerror}.
This fraction for $F$ 
is bounded by the fraction of $\bar{a} \in (A/\pp)^{np}$ 
such that $F$ and $\del F/\del t$ vanish mod $\pp$ at $\bar{a}$.
By Lemma~\ref{derivative}, $F$ and $\del F/\del t$
are relatively prime in $K[\dots,y_{ij},\dots]$,
so they define a subscheme of codimension at least $2$ in $\Aff^{np}_K$,
and the desired bound follows, 
just as in the last two sentences of Section~\ref{relativelyprimesection}.
\end{proof}

\section{Squarefree values of polynomials over other rings of functions}
\label{ringsoffunctions}

Let $S$ be a finite nonempty set of closed points of
a smooth, projective, geometrically integral curve $X$ over $\F_q$.
Define the affine curve $U:=X \setminus S$,
and let $A$ be the ring of regular functions on $U$.
Thus $A$ is the set of $S$-integers of the function field $K$ of $X$.
An element $a \in A$ is called squarefree if the ideal $(a)$ of $A$
is a product of distinct primes of $A$.
Let $\Div_S$ denote the set of effective divisors on $X$ with support
contained in $S$.
If $D \in \Div_S$, then the $\F_q$-subspace 
$L(D):=\{\, f \in K^*: (f)+D \ge 0 \,\} \cup\{0\}$ of $K$
is contained in $A$.
Let $\ell(D)=\dim_{\F_q} L(D)$.
Define the density $\mu(\calS)$ of a subset $\calS \subseteq A^n$ 
as the limit of $\#(\calS \cap \BOX)/\#\BOX$
as $\BOX$ runs through $L(D_1) \times \dots \times L(D_n)$,
where $D_i \in \Div_S$ and $\min_i \deg D_i \rightarrow \infty$.
The following theorem generalizes Theorem~1 in~\cite{ramsay1992};
we state it only for squarefree values, 
but as mentioned already in Section~\ref{theorems}
the lemmas used in its proof will yield the analogous result
for $k^{\operatorname{th}}$-power-free values.

\begin{theorem}[Squarefree values over rings of functions]
\label{ringsoffunctionstheorem}
With notation as above, 
let $f \in A[x_1,\dots,x_n]$ be a polynomial that is
squarefree as an element of $K[x_1,\dots,x_n]$.
Let
	$$\calS_f:=\{\, a\in A^n : \text{$f(a)$ is squarefree} \,\}.$$
For each nonzero prime $\pp \subseteq A$,
let $c_\pp$ be the number of $x \in (A/\pp^2)^n$
satisfying $f(x)=0$ in $A/\pp^2$.
Then $\mu(\calS_f)=\prod_{\pp} (1 - c_\pp/|\pp|^{2n})$,
where $|\pp|:=\#(A/\pp)$.
\end{theorem}

\begin{proof}[Sketch of proof of Theorem~\ref{ringsoffunctionstheorem}]
Most of the proof follows the proofs of Theorems
\ref{relativelyprime} and~\ref{Fqtsquarefree},
so we will comment only on the differences.

First we prove the analogue of Lemma~\ref{relativelyprimeerror}.
Enlarging $S$ only makes this analogue harder to prove,
since there are more boxes to consider,
while the primes $\pp$ involved are the same once $M$
exceeds $|\pp|$ for every new $\pp$ thrown into $S$.
Thus we may assume that $A$ is a principal ideal domain.
Let $g$ be the genus of $X$.
Given $\BOX = L(D_1) \times \dots \times L(D_n)$,
let $B_i = \deg D_i -(2g-2)$.
We may assume $\deg D_1 \le \dots \le \deg D_n$.
Replace the definition of $Q_s$ in the 
proof of Lemma~\ref{relativelyprimeerror} by
	$$Q_s:= \{\, a \in \BOX : \exists \pp \text{ such that 
	$B_s \le \deg \pp < B_{s+1}$ and $\pp \mid f(a),g(a)$} \,\}.$$
Replace the argument
``Since $|\pp|<B_{s+1}$,
the reductions modulo $\pp$
of the $(a_{s+1},\dots,a_n) \in \BOX(B_{s+1},\dots,B_n)$
are almost uniformly distributed in $(A/\pp)^{n-s}$''
by
``By Riemann-Roch, if $D$ is a divisor with $\deg D - \deg \pp > 2g-2$,
then $\ell(D)-\ell(D-\pp)=\deg \pp$;
hence the $\F_q$-linear reduction map 
	$$L(D_{s+1}) \times \dots \times L(D_n) \rightarrow (A/\pp)^{n-s}$$
is surjective for $\deg \pp < B_{s+1}$.''
Replace the estimate ``$|h(a_1,\dots,a_s)| = O(B_s^\eta)$''
by ``$\deg h(a_1,\dots,a_s) = O(\deg D_s) = O(B_s)$ where the 
implied constants depend on $h$ and the genus $g$, but not
on the $a_i$ or $D_i$''
to bound the number of primes $\pp$ with $\deg \pp \ge B_s$
dividing $h(a_1,\dots,a_s)$.
The rest of the proof of Lemma~\ref{relativelyprimeerror}
goes as before.

Next we prove the analogue of Lemma~\ref{Fqtsquarefreeerror}.
Choose any $t \in A-A^p$,
and note that Lemmas \ref{restrictionofscalars} and~\ref{derivative}
go through without change;
in fact, we can generalize Lemma~\ref{restrictionofscalars}
by replacing the form
$y_0^p+ty_1^p+\dots+t^{p-1}y_{p-1}^p$
by $t_1 y_1^p + t_2 y_2^p + \dots + t_m y^m$
for any $t_1,\dots,t_m \in A$ provided that $1$ and $t$ appear
among the $t_i$.

We claim that there exist $t_1,\dots,t_m \in A$
such that for every $D \in \Div_S$
and every $a \in L(D)$,
it is possible to write $a=t_1 a_1^p + \dots + t_m a_m^p$
with $a_i \in A$ such that $t_i a_i^p \in L(D)$ for each $i$.
More precisely, we claim that if
$E$ is the divisor $(p+1)(2g+5) \sum_{\pp \in S} \pp$,
then the choice $\{t_1,\dots,t_m\}:=L(E)$ works.
To prove this, it suffices to show, given $D \in \Div_S$
and $a \in L(D) - L(E)$, 
that one can adjust $a$ by an element of the form $t_i a_i^p$
in $L(D)$ to obtain an element in $L(D-\pp)$ for some $\pp$
in the support of $D$; then iterate.
If $a \in L(D) - L(E)$, then for some $\pp \in S$,
$a$ has a pole of order greater than $(p+1)(2g+5)$ at $\pp$.
If $n \ge 2g+4$,
Riemann-Roch shows that $\ell(n\pp)-\ell((n-1)\pp)=\deg \pp$;
in particular there exists a function in $A$ having a pole
of order $n$ with prescribed leading coefficient at $\pp$,
and no other poles.
Since every integer greater than $(p+1)(2g+5)$
is expressible as $i+pj$ with $2g+4 \le i < 2g+4+p$
and $j \ge 2g+4$,
we can find functions $t \in L(i\pp)$ and $\alpha \in L(j\pp)$
such that $t \alpha^p$ and $a$ have the same order of pole at $\pp$,
and the same leading coefficient.
Then $t \in L(E)$ and $t \alpha^p \in L(D)$, so this proves the claim.

By the previous paragraph, for each $D \in \Div_S$,
we have a surjection
\begin{align*}
	L(D_1) \times \cdots \times L(D_r) & \rightarrow L(D) \\
	(a_1,\dots,a_r) & \mapsto t_{i_1} a_1^p + \dots + t_{i_r} a_r^p
\end{align*}
for some subset $\{i_1,\dots,i_r\} \subseteq \{1,2,\dots,m\}$
and for some $D_i \in \Div_S$ depending on $D$.
Thus in proving the analogue of Lemma~\ref{Fqtsquarefreeerror},
we may reduce the result for $f$
to the result for each $F$ in the finite set of polynomials of the form
	$$F:=f \left( 
		\sum_{j \in S_1} t_j y_{1j}^p,
		\ldots,
		\sum_{j \in S_n} t_j y_{nj}^p \right) 
	\in K[\ldots,y_{ij},\ldots]$$
for all possible subsets $S_1,\dots,S_n \subseteq \{1,2,\dots,m\}$.

Choose $b \in A$
such that $b \frac{\del}{\del t}$ is a 
nonzero derivation $A \rightarrow A$.
Then the analogue of Lemma~\ref{Fqtsquarefreeerror} for an $F$ as above
follows from the analogue of Lemma~\ref{relativelyprimeerror}
applied to $F$ and $b \del F/\del t$.

To complete the proof of Theorem~\ref{ringsoffunctionstheorem},
it remains to prove $c_\pp = O(|\pp|^{2n-2})$,
to obtain convergence of the infinite product.
This follows as in the proof of Theorem~\ref{Fqtsquarefree}
at the end of Section~\ref{Fqtsection},
except that the fraction $c_\pp/|\pp|^{2n}$ for $f$
is now bounded by a {\em sum} of the analogous fractions 
for several different polynomials $F$.
Each of the latter fractions is bounded by $O(1/|\pp|^2)$,
and the number of polynomials $F$ is $O(1)$, independent of $\pp$,
so $c_\pp/|\pp|^{2n} = O(1/|\pp|^2)$, as desired.
\end{proof}

\section{The image of the values in $\Q^*/\Q^{*2}$}
\label{Qsquaredsection}

This section is devoted to the proof of Theorem~\ref{Qsquaredtheorem}.

\begin{lemma}
\label{thinset}
Suppose that $f(x) \in \Z[x]$ is squarefree as an element of $\Q[x]$,
and that $\deg f \ge 2$.
Fix $q \in \Q^*$ such that $q \not=1$.
Then the number of solutions $(m,n)$ to $f(m)= q f(n)$
satisfying $1 \le m, n \le B$
is $o(B)$ as $B \rightarrow \infty$.
\end{lemma}

\begin{proof}
Since $f(m)$ has constant sign for large positive $m$,
the result is trivial if $q<0$.
Therefore assume $q>0$.

Theorem~2 in Chapter~13 of~\cite{serremordellweil} implies that
the number of such solutions on each irreducible component of 
the curve $f(m)-qf(n)=0$ in the $(m,n)$-plane over $\Q$
is $O(B^{1/2} \log B)$ unless some component is a line.
If there is a line, it cannot be of the form $n=\alpha$
for any $\alpha \in \Q$,
so it would have an equation $m=\alpha n + \beta$
for some $\alpha,\beta \in \Q$.
Then $f(\alpha n +\beta) - q f(n)=0$ as polynomials in $\Q[n]$.
Equating leading coefficients shows that $\alpha\not=1$.
Choose $\gamma \in \Q$ so that $\alpha \gamma + \beta = \gamma$.
Then the polynomial $F(x):=f(x+\gamma)$ satisfies
$F(\alpha x) - q F(x)=0$.
Since $q>0$ and $q \not=1$, there is at most one integer $d$
such that $\alpha^d=q$.
Thus $F$ is a monomial of degree $d$,
and $f(n)=c(n-\gamma)^d$ for some $c \in \Q$.
Since $\deg f \ge 2$, this contradicts the assumption that
$f$ is squarefree.
\end{proof}

We now begin the proof of Theorem~\ref{Qsquaredtheorem}.
We easily reduce to the case where $f(x)=h(x)$,
that is, 
where the coefficients of $f$ have gcd~$1$,
and $f$ is squarefree in $\Z[x]$.
Also we may assume that the leading coefficient of $f$ is positive.
The $\deg f=0$ case is trivial.

\begin{proof}[Proof of Theorem~\ref{Qsquaredtheorem} in the case $\deg f=1$]
Write $f(x)=ax+b$.
Changing $b$ by a multiple of $a$
changes the sequence of values only in finitely many terms,
so we may assume $0<b \le a$.
Given $r,N \in \Z_{\ge 0}$,
let $S(r \bmod a, N)$ denote the set of positive squarefree integers $\le N$
that are congruent to $r$ modulo $a$.
Identify each element of $\Q^*/\Q^{*2}$ with a squarefree integer
representative.

We claim that the image of $\{f(1),f(2),\dots,f(B)\}$ in $\Q^*/\Q^{*2}$
is the disjoint union of $S(r \bmod a,\delta_r(aB+b))$
as $r$ ranges from $0$ to $a-1$.
Let $m$ be as in the definition of $\delta_r$,
when it exists.
If $s \in S(r \bmod a,\delta_r(aB+b))$,
then $\delta_r>0$, 
so the integer $m$ in the definition of $\delta_r$ exists;
then $m^2 s \equiv m^2 r \equiv b \pmod{a}$ and 
$m^2 s \le m^2 \delta_r(aB+b) = aB+b$,
so $m^2 s \in \{f(1),f(2),\dots,f(B)\}$,
and hence $s$ is in the image.
Conversely, suppose that the squarefree integer $s$
represents the image of $f(n)$ in $\Q^*/\Q^{*2}$
for some $n \in \{1,2,\dots,B\}$.
Thus $f(n) = \bar{m}^2 s$ for some $\bar{m}$.
Let $r \in [0,a)$ be such that $r \equiv s \pmod{a}$.
Then $\bar{m}^2 r \equiv \bar{m}^2 s = an+b \equiv b \pmod{a}$,
so the $m$ in the definition of $\delta_r$ exists,
and $m \le \bar{m}$.
Now $s = f(n)/\bar{m}^2 \le (aB+b)/m^2 = \delta_r (aB+b)$,
so $s \in S(r \bmod a,\delta_r(aB+b))$.
This proves the claim.

When $\gcd(a,r)=1$, the density of squarefree values of $ax+r$
equals
	$$\prod_{p \nmid a} \left( 1-p^{-2} \right) 
	= \frac{6}{\pi^2} \prod_{p \mid a} \left( 1-p^{-2} \right)^{-1},$$
so
	$$\# S(r \bmod a,N) = 
	\left( \frac{6}{\pi^2} \prod_{p \mid a} \left( 1-p^{-2} \right)^{-1} 
		+ o(1) \right) (N/a)$$
as $N \rightarrow \infty$.
The result follows upon setting $N=\delta_r (aB+b) \sim \delta_r aB$
for each $r \in [0,a)$ for which $\delta_r \not=0$
(such $r$ are necessarily prime to $a$),
and summing over $r$.
\end{proof}

\begin{proof}[Proof of Theorem~\ref{Qsquaredtheorem} in the case $\deg f \ge 2$]
Replacing $f(x)$ by $f(x+n)$ for some $n$,
we may assume that $0<f(1)<f(2)<\cdots$.
It suffices to show, given $\epsilon>0$,
that for sufficiently large $B$,
the image of $\{f(1),\dots,f(B)\}$ in $\Q^*/\Q^{*2}$
has size at least $(1-\epsilon)B$.
For a positive integer $n$, let $\fraks(n)$ denote the largest positive
integer $m$ such that $m^2 \mid n$.
Define
	$$S_m=\{\,n \in \{1,2,\dots,B\} : \fraks(f(n))=m \,\}.$$
By Lemma~\ref{squarefreeerror1}, 
the density of the set of integers $n$ such that $\fraks(f(n))$
is divisible by a prime $p>M$ tends to zero as $M \rightarrow \infty$.
Also, for each fixed prime $\ell$, the set
of integers $n$ such that $\ell^m \mid \fraks(f(n))$
tends to zero as $m \rightarrow \infty$,
because the number of $n \in \Z/\ell^{2m}$
such that $f(n)=0$ in $\Z/\ell^{2m}$ is $O(1)$ as $m \rightarrow \infty$,
since $f$ is squarefree.
Hence $\mu\left(\{\, n \in \Z: \fraks(f(n)) \ge M \,\}\right) \rightarrow 0$ 
as $M \rightarrow \infty$.
In particular, if $M$ is sufficiently large,
then $\#(S_1 \cup \dots \cup S_{M-1}) > (1-\epsilon/2)B$
for large $B$.

Now $f$ maps each $S_i$ injectively into $\Q^*/\Q^{*2}$.
For $1 \le i < j <M$,
the intersection of the images of $f(S_i)$ and $f(S_j)$
in $\Q^*/\Q^{*2}$ has size $o(B)$ as $B \rightarrow \infty$
by Lemma~\ref{thinset}.
Thus the image of $\{f(1),f(2),\dots,f(B)\}$
in $\Q^*/\Q^{*2}$ has size at least 
$(1-\epsilon/2)B - \binom{M-1}{2} o(B) > (1-\epsilon)B$,
if $B$ is sufficiently large relative to $M$.
\end{proof}

\section*{Acknowledgements}

I thank Johan de Jong and Ofer Gabber for independently suggesting
that I try to use ideas from~\cite{poonenbertini} 
to compute the density of squarefree values in the $\F_q[t]$ case.
I thank Brian Conrad for pointing out some errors in an earlier draft,
and the referee for helpful suggestions.
Finally, I thank Pamela Cutter, Andrew Granville, and Tom Tucker
for sharing their preprint with me.

\providecommand{\bysame}{\leavevmode\hbox to3em{\hrulefill}\thinspace}

\end{document}